\documentclass{amsart}

\usepackage{maple2e}
 
\DefineParaStyle{Maple Output}
\DefineParaStyle{Maple Plot}
\DefineCharStyle{2D Math}
\DefineCharStyle{2D Output}

\title{Zeta functions and regularized determinants on projective 
spaces}

\author{M. Spreafico}
\address{Dipartimento di Matematica ed Applicazioni, Universit\`a
Milano Bicocca, via Bicocca degli Arcimboldi 8, 20126 Milano, Italy} 
\date{}

\newtheorem{lem}{Lemma}

\newtheorem{defi}{Definition}

\newtheorem{prop}{Proposition}

\newcommand\e{{\rm e}}
\newcommand{\beq}{\begin{equation}}
\newcommand{\eeq}{\end{equation}}
\newcommand{\R}{{\rm Re}}

 \def\NN{\hbox{\sf I\kern-.13em\hbox{N}}}
 \def\HH{\hbox{\sf I\kern-.13em\hbox{H}}}
 \def\DD{\hbox{\sf I\kern-.13em\hbox{D}}}
 \def\RR{\hbox{\sf I\kern-.14em\hbox{R}}}
 \def\CC{\hbox{\sf I\kern-.44em\hbox{C}}}
 \def\ZZ{{\hbox{\sf Z\kern-.43emZ}}}
 \def\QQ{\hbox{\sf C\kern -.48emQ}}
 \def\Cc{\hbox{\sf C\kern -.47em {\raise .48ex \hbox{$\scriptscriptstyle |$}}
   \kern-.5em {\raise .48ex \hbox{$\scriptscriptstyle |$}} }}
 \def\Qq{\hbox{\sf Q\kern -.57em {\raise .48ex \hbox{$\scriptscriptstyle |$}}
   \kern-.55em {\raise .48ex \hbox{$\scriptscriptstyle |$}} }}
 \def\BB{\hbox{\sf ]\kern-.5em\hbox{[}}}

\begin{document}

\maketitle

\begin{abstract} A Hermite type formula is introduced and used to study the 
zeta function over the real and complex $n$-projective space. 
This approach allows to
compute the residua at the poles and the value at the origin as well
as the value of the derivative at the origin, that gives the
regularized determinant of the associated Laplacian operator.

\end{abstract}

\section{Introduction}
\label{s1}

Zeta functions on the sphere (and in general on a closed Riemannian
manifold) were first introduced by Minakshisundaram
and Peijel as extensions of the classical Riemann zeta function
\cite{MP}. An analytic definition, 
by a Mellin transform of the trace of the heat operator 
associated to the Laplacian in the standard metric, 
shows how the residua at the poles are given by the
coefficients in the asymptotic expansion of the trace of the heat
operator \cite{ABP}. Such an expansion can be obtained in a
large number of cases using global analysis \cite{See} \cite{Cal} \cite{Bru}
\cite{BS2} , 
and in particular the first coefficients in the case of 
the Laplacian on a Riemannian manifold can be computed using 
local invariants associated to the curvature tensor \cite {Gil}.
The constant term, that corresponds to the value of the zeta function at
the origin, is important in physics \cite{Wei} and in conformal
theory \cite{BO}, being
associated to the conformal anomaly. On the other hand, 
the derivative of the zeta function at the
origin gives the Atiyah regularized determinant of the Laplacian \cite{RS}. 
Early approaches to calculate these quantities give explicit results for the 2
sphere \cite{Wei} \cite{Dow}, 
while more recently an explicit formula for the residua has been  
obtained in \cite{CMB2} for the $n$-sphere, using a result that
allows to write a
Dirichlet series as a sum of classical Hurwitz zeta functions
\cite{CMB1}. This method fails to compute the derivative, but an
alternative one is provided by \cite{CQ}, 
where a factorization theorem for zeta
regularized products is introduced and explicit results for the 
low dimensional spheres are given.

The aim of this paper is to show how these zeta functions can be
treated very easily by classical methods using a 
Hermite type formula, in exactly the
same way as for the one dimensional case of the Riemann (Hurwitz) 
zeta function.
This approach allows us to deal not only with the zeta function on the spheres,
but also with the zeta function on the real and complex projective spaces. 
For all these cases, we localize the poles and give explicit
formulas for the residua. We show that the origin $s=0$ is a regular point, and
compute the value of the zeta function together with the value of its first
derivative at $s=0$.

The main feature of this approach is that it can be used more generally to deal
with the zeta function of any operator whose eigenvalues are explicitly known
with their multiplicity. There is work in progress where further cases are under
consideration.

\section{A Hermite type formula}
\label{s2}

Consider the function 
\[
z(s,a,b,c)=\sum_{n=1}^\infty \frac{P_d(cn)}{(cn+a)^s(cn+b)^s},
\]
of the complex variable $s$ for $\R (s)>\frac{d+1}{2}$, where $P_d$ is
a polynomial of degree $d$, and $a,b,c$ real constants with
$a,b>-1$, $c>0$.
Introduce the function of the complex variable $z$: 
\[
\phi(z,s,a,b,c)=\frac{P_d(z)}{(cz+a)^s(cz+b)^s},
\]
then we have the following: 
\begin{prop} For $\R(s)>\frac{d+1}{2}$,
\[
z(s,a,b,c)=\frac{1}{2}\phi(1,s,a,b,c)+\int_1^\infty \phi(z,s,a,b,c) dz+
I(s),
\]
where $I(s)$ is an integral function of $s$.
\label{l1}
\end{prop}
\noindent\underline{Proof} The infinite sum, is
actually a finite sum of infinite sums of the following type, that can be
treated by using the Plana theorem when $\R(s)>\frac{d+1}{2}$,
\[
\sum_{n=1}^\infty \psi_l(n,s,a,b,c)=
\sum_{n=1}^\infty \frac{c_l (cn)^l}{(cn+a)^s(cn+b)^s}=
\]
\[
=\frac{1}{2}\psi_l(1,s,a,b,c)+\int_1^\infty \psi_l(z,s,a,b,c) dz+
\]
\[
-2c_l\int_0^\infty
(c^2+y^2)^{l/2}[(a+c)^2+y^2]^{-s/2}[(b+c)^2+y^2]^{-s/2}\times 
\]
\[
\times
\sin \left[l\arctan\frac{y}{c}-s\left(
\arctan\frac{y}{a+c}+\arctan\frac{y}{b+c}\right)\right]
\frac{dy}{e^{2\pi y}-1},
\]

Indeed, the last integral converges for all values of $s$.
Furthermore, the last integral converges 
uniformly by standard estimates on $\arctan$, so the given formula
defines an analytic function of $s$.
$\BB$

From proposition \ref{l1}, it is clear that all the poles and relative
residua come from the first two terms. This expression can also be used
to compute the value of $z$ and its derivative at $s=0$.

\section{Formulae for the residua and the value at the origin}
\label{s3}

Spheres and real projective spaces can be treated together as follows, complex
projective spaces will be considered afterwards.

The polynomial $P_d$, with $c=1,2$ for 
the sphere $S^k$ and the projective space $\RR P^k$ respectively ($k>1$), is:
\[
P_{k-1}(x)=Q_k(x)= \frac{2x+k-1}{(k-1)!}\prod_{i=1}^{k-2} (x+i).
\]

It is now convenient to distinguish odd and even cases. We introduce
the numbers $b_{k,l}$, by

\begin{defi} For $k=1,2, 3, \dots$, 
the numbers $b_{k,l}$ are defined by the equations
\[
\prod_{i=1}^{2h-2} (x+i)=\sum_{l=0}^{h-1} b_{2h,l} [x^2+(2h-1)x]^l,
\]
\[
(x+h)\prod_{i=1}^{2h-1} (x+i)=\sum_{l=0}^{h} b_{2h+1,l} [x^2+2hx]^l.
\]
\label{d1}
\end{defi}

In particular: $b_{2h,h-1}=b_{2h+1,h}=1$, $b_{2h,0}=(2h-2)!$,
$b_{2h+1,0}=h(2h-1)!$. 

\vskip .2in

This allows us to write the zeta function as a sum of some standard
ones. Let  
\[
z_{k}(s,c)=\sum_{n=1}^\infty \frac{Q_k(cn)}{[cn(cn+k-1)]^s},
\]
be the zeta function in dimension $k$; then
\[
\zeta(s,S^k)=z_k(s,1),
\]
\[
\zeta(s,\RR P^k)=z_k(s,2).
\]

Now, consider the two functions
\[
z_{even}(s,a,b,c)=\sum_{n=1}^\infty \frac{2(cn+a)+b-a}{[(cn+a)(cn+b)]^s},
\]
\[
z_{odd}(s,a,b,c)=\sum_{n=1}^\infty \frac{1}{[(cn+a)(cn+b)]^s},
\]
then:
\begin{lem} For $h=1, 2, \dots$, 
\[
z_{2h}(s,c)=\frac{1}{(2h-1)!}\sum_{l=0}^{h-1} b_{2h,l}
z_{even}(s-l,0,2h-1,c), 
\]
\[
z_{2h+1}(s,c)=\frac{2}{(2h)!}\sum_{l=0}^{h-1} b_{2h+1,l}
z_{odd}(s-l,0,2h,c). 
\]
\label{l3}
\end{lem}

To write simpler formulae, let us also introduce the function
\[
f(s,l,a,b,c)=\int_0^\infty [(a+c)^2+y^2]^{(l-s)/2}
[(b+c)^2+y^2]^{-s/2}\times \]
\[
\times
\sin
\left[(l-s)\arctan\frac{y}{a+c}-s\arctan\frac{y}{b+c}\right]\frac{dy}{e^{2\pi
y}-1}.
\]

By expanding the trigonometric functions and recalling the standard
integral representation of the Bernoulli numbers $B_n$, we
get\footnote{Observe that $f(0,0,a,b,c)=0$.}
\[
f(-n, l, a,b,c)=
\]
\[
=\frac{1}{4}\sum_{i=1}^{E\left(\frac{n+l+1}{2}\right)} 
\sum_{j=0}^{E\left(\frac{n}{2}\right)} 
\left(\begin{array}{c}n+l\\2i-1\end{array}\right)
\left(\begin{array}{c}n\\2j\end{array}\right)
(a+c)^{n+l+1-2i} (b+c)^{n-2j} \frac{B_{2(i+j)}}{i+j}+
\]
\[
+\frac{1}{4}\sum_{i=1}^{E\left(\frac{n+1}{2}\right)} 
\sum_{j=0}^{E\left(\frac{n+l}{2}\right)} 
\left(\begin{array}{c}n\\2i-1\end{array}\right)
\left(\begin{array}{c}n+l\\2j\end{array}\right)
(a+c)^{n+l-2j} (b+c)^{n+1-2i} \frac{B_{2(i+j)}}{i+j},
\]
where $E(q)$ denotes the integer part of the rational number $q$.

We can now state the main properties of the functions $z_{even/odd}$.

\begin{lem} The function $z_{even}(s,a,b,c)$ has a simple pole at
$s=1$ with residuum $\frac{1}{c}$, while for $m=0,1,2,\dots$,
\[
z_{even}(-m,a,b,c)=\frac{1}{2}[2(c+a)+b-a](c+a)^m(c+b)^m+
\]
\[
-\frac{(c+a)^{m+1}(c+b)^{m+1}}{(m+1)c}-4f(-m,1,a,b,c)-2(b-a)f(-m,0,a,b,c).
\]
\label{l4}
\end{lem}
\noindent\underline{Proof} 
Proceeding as in proposition \ref{l1}, we get 
\[
z_{even}(s,a,b,c)=\frac{1}{2}\frac{2(c+a)+b-a}{(c+a)^s(c+b)^s}+\frac{1}{c}
\frac{1}{(c+a)^{s-1}(c+b)^{s-1}}\frac{1}{s-1}+
\]
\[
-4\int_0^\infty [(a+c)^2+y^2]^{1/2-s/2}[(b+c)^2+y^2]^{-s/2}\times
\]
\[
\times
\sin\left[(1-s)\arctan\frac{y}{a+c}-s\arctan\frac{y}{b+c}\right]
\frac{dy}{e^{2\pi y}-1}+
\]
\[
+2(b-a)\int_0^\infty [(a+c)^2+y^2]^{-s/2}[(b+c)^2+y^2]^{-s/2}\times
\]
\[
\times
\sin\left[s\left(\arctan\frac{y}{a+c}+\arctan\frac{y}{b+c}\right)\right]
\frac{dy}{e^{2\pi y}-1};
\]
then the poles are given by the second term, while the values at the non
positive integers can be computed by using the formulae introduced above.
$\BB$

\begin{lem} The function $z_{odd}(s,0,b,c)$ has simple poles at
$s=\frac{1}{2}-m$, $m=0, 1, 2, \dots$, with residua 
\[
{\rm Res}_1\left(z_{odd}(s,0,b,c), s=\frac{1}{2}-m\right)
=\frac{(-1)^m}{2^{m+1}}
\frac{(2m-1)!!}{m!}\frac{1}{c}\left(\frac{b}{2}\right)^{2m}.
\]

For $m=0,1,2,\dots$, 
\[
z_{odd}(-m, 0,b,c)=\frac{1}{2}c^m
(c+b)^m+\frac{(-1)^{m+1}}{m+1}\frac{2^m}{c}\left(\frac{b}{2}\right)^{2m-1}
\frac{\Gamma(m+1)}{(2m+1)!!}+
\]
\[
-\frac{(bc)^m}{m+1}\frac{\Gamma(2m+1)}{\Gamma(-m)\Gamma(m+1)}
\sum_{i=0}^m
\frac{\Gamma(i-m)\Gamma(i+m+1)}{i!\Gamma(2m+i+1)}\left(\frac{c}{b}\right)^i
-2f(-m,0,0,b,c).
\]
\label{l5}
\end{lem} 
\noindent\underline{Proof}
For the odd case, difficulties arise when dealing with the first
integral. If we restrict ourselves to the interesting case of $a=0$, we get
\[
z_{odd}(s,0,b,c)=\frac{1}{2}\frac{1}{c^s(c+b)^s}
+c^{-s}\int_1^\infty x^{-s}(cx+b)^{-s} dx+
\]
\[
+2\int_0^\infty [c^2+y^2]^{-s/2}[(b+c)^2+y^2]^{-s/2}
\sin\left[s\left(\arctan\frac{y}{c}+\arctan\frac{y}{b+c}\right)\right]
\frac{dy}{e^{2\pi y}-1},
\]
and we have the following two possible ways of treating it:
\[
c^{-s}\int_1^\infty x^{-s}(cx+b)^{-s} dx
=c^{-s}(c+b)^{-s} \frac{1}{2s-1} F\left(s,1;2s;\frac{b}{c+b}\right)=
\]
when $\R(s)>\frac{1}{2}$ and where $F$ is the hyper geometric function, and
\[
=\frac{1}{\sqrt\pi}\frac{1}{c}\left(\frac{b}{2}\right)^{1-2s}
\frac{\Gamma(1-s) \Gamma\left(s+\frac{1}{2}\right)}{2s-1} -
c^{-s}\int_0^1 x^{-s}(cx+b)^{-s} dx,
\]
when $\frac{1}{2}<\R(s)<1$, where the last integral is actually
convergent if $\R(s)<1$, and can be expressed in term of
a hyper geometric function:
\[
c^{-s}\int_0^1 x^{-s}(cx+b)^{-s} dx=\frac{(bc)^{-s}}{1-s}
F\left(s,1-s;2-s;-\frac{c}{b}\right). 
\]

In particular, we use the first
representation to get the analytical continuation in the negative half plane
and to compute the residua at the poles, and the second one
to compute the value at $s=-m$, $m=0,1,2,\dots$, 
and the derivative at $s=0$ (see section \ref{s4})
(Both expressions are good to compute the residua).
$\BB$

From the previous lemmas we immediately get the following:

\begin{prop} The function $z_{2h}$ has simple poles at $s=n$, for
$n=1,2,3,\dots, h$, with residua:
\[
{\rm Res}_1(z_{2h}(s,c), s=n)=\frac{1}{c}\frac{b_{2h, n-1}}{(2h-1)!}.
\]

The function $z_{2h+1}$ has simple poles at $s=\frac{1}{2}+h-m$,
$m=0,1,2,\dots$, with residua
\[
{\rm Res}_1\left(z_{2h+1}(s,c), s=\frac{1}{2}+h-m\right)=
\]
\[
=\frac{1}{c}\frac{2}{(2h)!}\sum_{l=0}^{{\rm min}(h,m)}
\frac{(-1)^{m-l}b_{2h+1,l}}{2^{m-l+1}} 
\frac{(2(m-l)-1)!!}{(m-l)!} \left(\frac{b}{2}\right)^{2(m-l)}.
\]
\label{p2}
\end{prop}

These should be compared with theorem 2 of \cite{CMB2}, where we note that a
factor $\frac{1}{2}$ is missing.

For what concerns the value of $z_k(0,c)$, this can be easily computed
(using any mathematical software) from the formulae in lemmas \ref{l3},
\ref{l4} and \ref{l5}. Explicit results for $n=2,3$ and $4$ are given
in section \ref{s4}.

Next, consider the complex projective spaces $\CC P^k$, $k>1$. The polynomial
$P_d$ is
\[
P_{k-1}(x)=Q_k(x)=\frac{k(2x+k)}{[(k-1)!]^2} \prod_{i=1}^{k-1} (x+i)^2,
\]
and 
\[
\zeta(s,\CC P^k)=\frac{1}{4^s}\sum_{n=1}^\infty \frac{Q_k(x)}{n^s(n+k)^s}.
\]

If we introduce the numbers $a_{k,l}$ by

\begin{defi} For $k=1,2, 3, \dots$, 
the numbers $a_{k,l}$ are defined by the equations
\[
\prod_{i=1}^{k-1} (x+i)^2=\sum_{l=0}^{k-1} a_{k,l} (x^2+kx)^l,
\]
\label{d2}
\end{defi}
(in particular: $a_{k,k-1}=1$, $a_{k,0}=[(k-1)!]^2$); then we can write

\begin{lem}
\label{l6}
For $k=1,2,3,\dots$, 
\[
\zeta(s,\CC P^k)=\frac{1}{4^s}\frac{k}{[(k-1)!]^2} \sum_{l=0}^{k-1} a_{k,l}
z_{even}(s-l,0,k,1).
\]
\end{lem}

From this and lemma \ref{l4}, we immediately get 

\begin{prop}
\label{p3}
The function $\zeta(s,\CC P^k)$ has simple poles at $s=n$, for $n=1,2,\dots,
k$, with residua
\[
{\rm Res}_1(\zeta(s, \CC P^k), s=n)=\frac{1}{4^n} \frac{k}{[(k-1)!]^2}
a_{k,n-1}.
\]
\end{prop}

\section{Determinant and low dimensional cases}
\label{s4}

In this last section we show how to use proposition \ref{l1} to compute
the regularized determinant by considering the case of $k=3$. Next, we
give explicit results for dimensions 2,3 and 4. Even if this is beyond the
purposes of these notes, we observe that such results show that a unique
formula for the determinant for a general $k$ is likely to exist, and 
could be determined using the approach outlined here.

We have $k=3$, $h=1$, $a=0$ and $b=2$.
Writing $z_3(s,c)$ as in proposition \ref{l1}, we decompose the second
term as in lemma \ref{l3}, but the third one by developing the
numerator. We get
\[
z_3(s,c)
=\sum_{n=1}^\infty \frac{(cn+1)^2}{[cn(cn+2)]^s}=
\]
\[
=\frac{(c+1)^2}{2[c(c+2)]^s} +\int_1^\infty [cx(cx+2)]^{-s}dx+
\int_1^\infty [cx(cx+2)]^{1-s}dx+
\]
\[
-2\int_0^\infty(c^2+y^2)^{-s/2}[(c+2)^2+y^2]^{-s/2}\sin[-s(\theta+\phi)]
\frac{dy}{e^{2\pi y}-1}+
\]
\[
-4\int_0^\infty(c^2+y^2)^{1/2-s/2}[(c+2)^2+y^2]^{-s/2}\sin[(1-s)\theta-s\phi]
\frac{dy}{e^{2\pi y}-1}+
\]
\[
-2\int_0^\infty(c^2+y^2)^{1-s/2}[(c+2)^2+y^2]^{-s/2}\sin[(2-s)\theta-s\phi]
\frac{dy}{e^{2\pi y}-1},
\]
where $\theta=\arctan \frac{y}{c}$, $\phi=\arctan \frac{y}{c+1}$.
The derivative at $s=0$ of the first term is immediate, 
\[
\left.
\frac{d}{ds}\frac{(c+1)^2}{2[c(c+2)]^s}\right|_{s=0}=-(c+1)\log\sqrt{c(c+1)},
\]
while that
of the second one can be computed with little effort by using the
representation introduced in the proof of lemma \ref{l5}. We get
\[
-\frac{8}{9c}-\frac{4}{3}-\frac{2c}{3}-\frac{2c^2}{9}
+\left(\frac{1}{3c}+1+c+\frac{c^2}{3}\right)\log c+
\left(\frac{2}{3c}+1+c+\frac{c^2}{3}\right)\log
\left(1+\frac{2}{c}\right). 
\]

For the last three integrals we get
\[
2\int_0^\infty (\theta+\phi)\frac{dy}{e^{2\pi y}-1}+
\]
\[
+2\int_0^\infty\sqrt{c^2+y^2}\{\log(c^2+y^2)+\log[(c+2)^2+y^2]\}\sin\theta
\frac{dy}{e^{2\pi y}-1}+
\]
\[ 
+4\int_0^\infty\sqrt{c^2+y^2}(\theta+\phi) \cos\theta 
\frac{dy}{e^{2\pi y}-1} +
\]
\[
+\int_0^\infty(c^2+y^2)\{\log(c^2+y^2)+\log[(c+2)^2+y^2]\}\sin 2\theta
\frac{dy}{e^{2\pi y}-1} +
\]
\[
+2\int_0^\infty\sqrt{c^2+y^2}(\theta+\phi) \cos 2\theta 
\frac{dy}{e^{2\pi y}-1} =
\]
that can be expressed in terms of derivatives of the Hurwitz
zeta function,
\[
=\zeta_H'(-2,c)+\zeta_H'(-2,c+2)+2\zeta_H'(-1,c)-2\zeta_H'(-1,c+2)+
\]
\[
+(3-2c)(\zeta_H'(0,c)+\zeta_H'(0,c+2))+\frac{44}{9}+\frac{4}{3}c
-\frac{10}{3}c^2+\frac{2}{9}c^3+
\]
\[
+\left(\frac{3}{2}-3c+\frac{3}{2}c^2-\frac{1}{3}c^3\right)\log c
+\left(-\frac{19}{6}+\frac{3}{2}c^2+c-\frac{1}{3}c^3\right)\log (c+2).
\]

Collecting and simplifying, we get the results shown below, together
with the other low dimensional cases.

\[
\zeta(0,S^2)=-\frac{2}{3}=-0.\bar 6,
\]
\[
\zeta'(0,S^2)=4\zeta'_R(-1)-\frac{1}{2}=-1.162 \dots;
\]
\[
\zeta(0,\RR P^2)=-\frac{11}{3}=-3.667 \dots,
\]
\[
\zeta'(0,\RR P^2)=4\zeta'_R(-1)-3\log 3+\frac{5}{2}=-1.459 \dots;
\]
\[
\zeta(0,S^3)=-1,
\]
\[
\zeta'(0,S^3)=2\zeta'_R(-2)+2\zeta'_R(0)+\log 2=-1.206 \dots;
\]
\[
\zeta(0,\RR P^3)=-\frac{1}{2}=-0.5,
\]
\[
\zeta'(0,\RR P^3)=2\zeta'_R(-2)-2\zeta'_R(0)-\frac{13}{6}\log 2+2\log
3-8=-5.527 \dots;
\]
\[
\zeta(0,S^4)=-\frac{61}{90}=-0.6\bar 7,
\]
\[
\zeta'(0,S^4)=\frac{2}{3}\zeta'_R(-3)+\frac{13}{3}\zeta'_R(-1)+\log 3 
-\frac{15}{16}=-0.5516 \dots,
\]
\[
\zeta(0,\RR P^4)=\frac{7}{45}=0.1\bar 5,
\]
\[
\zeta'(0,\RR P^4)=\frac{2}{3}\zeta'_R(-3)+\frac{13}{3}\zeta'_R(-1)
+\frac{13}{6}\log 2+\log 3 -\frac{35}{16}\log 5+\frac{45}{16}=-4.684 \dots,
\]
\[
\zeta(0,\CC P^2)=-\frac{89}{30}=-2.9\bar 6,
\]
\[
\zeta'(0,\CC P^2)=8 \zeta'_R(-3)+24 \zeta'_R(-1)+\frac{149}{15}\log 2 -4 \log 3
-\frac{203}{12}=-18.353\dots.
\]

The values for $\zeta'(0,S^k)$ agree with the ones provided by \cite{CQ},
section 4; the value for $\zeta(0,S^2)$ confirms the one originally given
by \cite{Wei}, against that of \cite{Dow}(see also section \ref{s5} for further
remarks); the values for $\zeta(0,S^3)$
agrees with the one computed from the formula in the corollary of section 3
of \cite{CMB2} (see also section 3 of \cite{CMB1}). Numerical
computations were done by using Maple.

\section{Remarks}
\label{s5}

To conclude some remarks are in order. The first one concerns the one
dimensional case, where we get the following relations
\[
\zeta(s,S^1)=2\zeta_R(2s),
\]
\[
\zeta(s,\RR P^1)=2^{-2s} \zeta(s,S^1),
\]
\[
\zeta(s,\CC P^1)=2^{-2s}\zeta(s,S^2).
\]

The first two relations are well known, while the third one can be easily read
out from the results in the sections above. We can then complete the table in
section \ref{s4}:
\[
\begin{array}{cc}\zeta(0,S^1)=-1,&\zeta'(0,S^1)=4\zeta'_R(0)=-3.676\dots;\\
\zeta(0,\RR P^1)=-1,&\zeta'(0,\RR P^1)=4\zeta'_R(0)+2\log 2=-2.29\dots;\\
\zeta(0,\CC P^1)=-\frac{2}{3}=-0.\bar 6,&\zeta'(0,\CC
P^1)=4\zeta'_R(0)+\frac{4}{3}\log 2-\frac{1}{2}=-3.252\dots.
\end{array}
\]

Also notice that, while the value of the zeta function at $s=0$ seems to depend
only on the topology, the value of its derivative, and hence the
regularized determinant, does not. For what concerns the first statement, just 
recall 
that, for any closed Riemannian manifold $M$ of dimension $m$,
\[
\zeta(0,M)=a_m(M)-{\rm dimker} \Delta_M,
\]
where $\Delta_M$ is the Laplacian operator in the standard metric, and $a_m(M)$
is the coefficient of the constant term in the asymptotic expansion of the heat
operator $\e^{-t\Delta_M}$ (see for example \cite{Ros}). In particular, compare
with the value computed above for $\zeta(0,S^2)$, where 
$a_2(S^2)=\frac{1}{24\pi}\int_{S^2} R_{S^2}(x) dvol(x)=\frac{1}{3}$.

Our final remark, concerns the possibility of using the approach introduced
to treat the case of the Laplacian coupled with a constant potential.
This is an important tool to face the problem of a generic potential (see
\cite{LS} for the one dimensional case). For the sake of simplicity, we restrict
ourselves to the case of the 2-sphere. The eigenvalues are then
$\lambda_n=n(n+1)+q^2$, and the zeta function is
\[
\zeta(s,\Delta_{S^2}+q^2)=\sum_{n=1}^\infty \frac{2n+1}{[n(n+1)+q^2]^s}.
\]

By expanding the power of the binomial (for finite $q$), this becomes
\[
\zeta(s,S^2)-s\zeta(s+1,S^2)q^2+\sum_{k=2}^\infty \left(\begin{array}{c}-s\\
k\end{array}\right)  \zeta(s+k,S^2)q^{2k}.
\]

We can compute
\[
\zeta(0,\Delta_{S^2}+q^2)=\zeta(0,S^2)-q^2,
\]
\[
\zeta'(0,\Delta_{S^2}+q^2)=\zeta'(0,S^2)+(1-2\gamma)q^2-\log G(q),
\]
where $\gamma=-\psi(1)$, and $G(q)$ is the integral function of order 2 
and zeros
$\pm i\sqrt{n(n+1)}$, of multiplicity $2n+1$, defined by the canonical product of
genus 2:
\[
G(q)=\prod_{n=1}^\infty \left[1+\frac{q^2}{n(n+1)}\right]^{2n+1}
\e^{-\frac{2n+1}{n(n+1)} q^2}.
\]

The main difference with respect to the one dimensional case is due to the
slower rate of increasing of the eigenvalues; this reflects in the appearance of
a further contribution coming from the singularity of $\zeta(s,S^2)$ at $s=1$,
thorough the term $s\zeta(s+1,S^2)$. In the following figures are plotted the
functions $G(q)$ and $\zeta'(0,\Delta_{S^2}+q^2)$, where the infinite product is
approximated by the product on the first 1000 terms.

\begin{maplegroup}
\mapleresult
\begin{center}
\mapleplot{G03.eps}
\end{center}
\end{maplegroup}
\centerline{Fig. 1, $G(q)$}

\mapleresult
\begin{maplegroup}
\mapleresult
\begin{center}
\mapleplot{det01.eps}
\end{center}
\end{maplegroup}

\centerline{Fig. 2, $\zeta'(0,\Delta_{S^2}+q^2)$}


\begin{thebibliography}{99}



\bibitem{ABP} M. Atiyah, R. Bott and V.K. Patodi, {\em On the Heat
Equation and the Index Theorem}, Inventiones math. 19 (1973) 279-330;

\bibitem{BO} Branson and Oersted, {\em Conformal geometry and local
invariants}, Diff. Geom. Appl. 1 (1991) 279-308;

\bibitem{Bru} J. Bruening, {\em Heat equation asymptotics for singular
Sturm-Liouville operators}, Math. Ann. 268 (1984) 173-196;


\bibitem{BS2} J. Bruening and R. Seeley, {\em The resolvent expansion for
second order regular singular operators}, J. of Funct. An. 73 (1988)
369-415; 

\bibitem{Cal} C. Callias, {\em The heat equation with singular
coefficients}, Comm. Math. Phys. 88 (1983) 357-385;

\bibitem{CMB1} E. Carletti and G. Monti Bragadin, {\em On
Dirichlet series associated with polynomials}, Proc. Am. Math.
Soc. 121 (1994) 33-37;

\bibitem{CMB2} E. Carletti and G. Monti Bragadin, {\em On
Minakshisundaram-Pleijel zeta functions on spheres}, Proc. Am. Math.
Soc. 122 (1994) 993-1001;

\bibitem{CQ} J. Choi and J.R. Quine, {\em Zeta regularized products and
functional determinants on spheres}, Rocky Mount. Jour. Math. 26 (1996)
719-729;

\bibitem{Dow} J.S. Dowker, {\em Effective actions in spherical
domains}, Comm. Math. Phys. 162 (1994) 633-647;

\bibitem{Gil} P.B. Gilkey, {\em Invariance theorems, the heat equation,
and the Atiyah-Singer index theorem}, Studies in Adv. Math. CRC Press, 1995;



\bibitem{LS} S. Levit and U. Smilansky, {\em A theorem on infinite products of
eigenvalues of Sturm-Liouville type operators}, Proc. AMS 65 (1977) 299-302;


\bibitem{MP} S. Minakshisundaram and A. Pleijel, {\em Some properties
of the eigenfunctions of the Laplace operator on Riemannian manifolds},
Canad. J. Math. 1 (1949) 242-256;

\bibitem{RS} D.B. Ray and I.M. Singer, {R-torsion and the Laplacian on
Riemannian manifolds}, Adv. Math. 7 (1974) 145-210;

\bibitem{Ros} S. Rosenberg, {\em The Laplacian on a Riemannian Manifold}, London
Math. Soc. ST31 1997;


\bibitem{Wei} W.I. Weisberger, {\em Conformal invariants for determinants
of Laplacians on Riemannian surfaces}, Comm. Math. Phys. 112 (1987) 633-638;


\bibitem{See} R.T. Seeley, {\em Complex powers of an elliptic operator},
Singular Integrals (Proc. Sympos. Pure Math. Chicago) 188-307, 
Amer. Math. Soc. (1967);





\end{thebibliography}
\end{document}